\documentclass{article}

\usepackage{personal-modified-frozen}

\theoremstyle{definition}
\newtheorem*{remark}{Remark}

\newminted[code]{idris}{
    gobble=8,
    breaklines,
    frame=lines,
    framesep=1em,
    rulecolor=\color{gray},
    linenos
}

\newmintinline[idris]{idris}{}

\newcommand \strong [1] {\textbf{#1}}

\newcommand \Type {\mathrm{Type}}

\newcommand \Judgement {\mathcal{J}}
\DeclareMathOperator  \Inference {\mathsf{Inference}}

\DeclareMathOperator    \Ruleset {Ruleset}
\DeclareMathOperator  \Deduction {\mathsf{Deduction}}

\newcommand*{\true}{\relax\ifnum\lastnodetype>0 \mskip\medmuskip\fi\mathsf{true}}
\newcommand*{\prop}{\relax\ifnum\lastnodetype>0 \mskip\medmuskip\fi\mathsf{prop}}
\newcommand*{\type}{\relax\ifnum\lastnodetype>0 \mskip\medmuskip\fi\mathsf{type}}

\DeclareMathOperator     \Set {Set}

\usepackage[all]{xy}

\newcommand \universe {\mathcal{U}}

\newcommand \Sig {\mathrm{Sig}}
\DeclareMathOperator \Sen {Sen}
\newcommand \Ded {\mathcal{D}}
\newcommand \Prop {\mathrm{P}}

\newcommand \connective {\mathcal{C}}
\DeclareMathOperator \arity {{\#}}
\DeclareMathOperator \rules {\mathcal{R}}
\DeclareMathOperator \headers {\mathcal{H}}
\newcommand \var {\mathcal{V}}
\newcommand \ty {\mathsf{ty}}
\newcommand \tm {\mathsf{tm}}

\DeclareMathOperator \application {\varepsilon}
\DeclareMathOperator \InL {\mathsf{In}_\mathrm{L}}
\DeclareMathOperator \InR {\mathsf{In}_\mathrm{R}}
\newcommand \caseof [5] {\mathsf{case}\;#1\;\mathsf{of}\;\InL#2\mapsto#3;\;\InR#4\mapsto#5}
\newcommand \Left {\mathsf{Left}}
\newcommand \Right {\mathsf{Right}}
\newcommand \thematicbreak {\begin{center}* * *\end{center}}
\DeclareMathOperator \term {\mathsf{term}}
\DeclareMathOperator \deduction {\mathsf{deduction}}

\title{
    Formalizing the \\
    Curry-Howard Correspondence}
\author{
    Juan F. Meleiro \\
    \texttt{juan.meleiro@me.com}
    \and Hugo L. Mariano \\
    \texttt{hugomar@ime.usp.br}
}
\date{2019}

\begin{document}

    \maketitle

    \begin{abstract}
        The Curry-Howard Correspondence has a long history, and still is a topic
        of active research. Though there are extensive investigations into the
        subject, there doesn't seem to be a definitive formulation of this
        result in the level of generality that it deserves. In the current work,
        we introduce the formalism of p-institutions that could unify previous
        aproaches. We restate the tradicional correspondence between typed
        $\lambda$-calculi and propositional logics inside this formalism, and
        indicate possible directions in which it could foster new and more
        structured generalizations.

        Furthermore, we indicate part of a formalization of the subject in the
        programming-language Idris, as a demonstration of how such
        theorem-proving enviroments could serve mathematical research.

        \smallskip
        
        \textbf{Keywords}. Curry-Howard Correspondence, p-Institutions, Proof
        Theory.
    
    \end{abstract}

    \tableofcontents

    \section{Some things to note}

    This work is the conclusion of two years of research, the first of them
    informal, and the second regularly enrolled in the course MAT0148
    \textit{Introdução ao Trabalho Científico}. As a final product, it's not a
    complete picture of the process that gave its origin. During that time, many
    subjects were discussed that didn't make their way into this pages -- be it
    for time, be it for a simple thematic inadequacy. Some of these
    would-be-results are mentioned along the way, and some have lost themselves
    to time or memory -- maybe to br remembered later on.

    Our investigations had simple motivations: clarity. The Curry-Howard
    Correspondence, our object of study, traces its origin to the early
    twentieth-century, with the invention of the $\lambda$-calculus, and is an
    active research subject even today. Being a correspondence between two
    logical systems, there is always the pursuit, inherent to mathematical
    research, to \emph{extend} its limits to more and more general classes.
    There are also situations, such as in the foudational tome \citet{HoTT2013},
    in whose foundational theoretical and methodological core is ingrained the
    Curry-Howard Correspondence. This idea is also known as the paradigm of
    \emph{propositions-as-types}.

    It's enough to say that the Curry-Howard Correspondence is as much inspiring
    as it is mysterious. And it's not for nothing that it doesn't have an
    stablished formalization. There are diverse formalisms and extensions built
    atop it, as the Lambda Cube \citep[see][]{BarendregtGeneralized1991}, or
    the formalism of \citet{LambekCategorical1988}. Each of them, however,
    has its compromises and imparsimonies\footnotemark.

    \footnotetext{Lambek, for example, restricts himself to some particular
    kinds of logics, while the Lambda Cube is fundamentally limited as a
    mechanism for extending Curry-Howard, as it deals only with
    \emph{subsystems} of the Calculus of Inductive Constructions.}

    So we decided to carve a sliver into this problem ourselves, and seek our
    own understanding of what is this result. And of course, what we found is a
    reflection of our own methodological aproach -- a personal expression.
    Nonetheless, we oriented ourselves with the goal of distiling the essence of
    the problem.

    \subsection{What we are talking about}

    The Curry-Howard Correspondence is an observation on the relation between
    two deductive systems: the intuitionistic propositional logic and the simply
    typed lambda calculus.

    The propositional logic is familiar to any mathematician. It's the logic
    that deals, evidently, with propositions -- these formal objects that
    represent statements and ``it is raining \emph{and} it is cold'', or
    ``\emph{if} it is raining \emph{them} i'll take my umbrella''. They are
    representad symbolically by abstracting details off the non-propositional
    parts with variables $A$, $B$, $P$, $Q$, \etc.

    From an algebraic point of view, such logics are algebras with operations
    such as the conjunction, disjunction, implication, \etc.

    From a logical point of view, there are multiple formalisms for such
    systems, but we'll use the one in which proofs look like trees

    \[
        \inferrule {
            \inferrule {
                \inferrule {
                    \inferrule {
                        \inferrule { } {
                            \vdash A \prop
                        }
                    } {
                        A \vdash A \true
                    } \\
                    \inferrule {
                        \inferrule { } {
                            \vdash A \prop
                        } \\
                        \inferrule { } {
                            \vdash B \prop
                        }
                    } {
                        A \vdash B \prop
                    }
                } {
                    A, B \vdash A \true
                }
            } {
                A \vdash B \rightarrow A \true
            }
        } {
            \vdash A \rightarrow (B \rightarrow A) \true
        }
    \]

    The $\lambda$-calculus, on the other hand, has its origin in Computer
    Science, offering (in its untyped version) an alternative to Turing
    Machines.

    The terms of the $\lambda$-calculus follow a specific syntax.

    \[
        \term = \var \mathrel| (\term \term) \mathrel| \lambda \var \term
    \]

    which represent, respectively, variables, the application of functions, and
    the construction of functions. For example, the \emph{constant} function has
    the form

    \[
        \lambda x \lambda y \mathrel. x
    \]

    which takes and argument $x$ and returns a function that, upon taking any
    other argument $y$, always returns $x$.

    The crucial observation, here, is that we can make this theory more
    \emph{well-behaved} by attributing to this terms \emph{types} -- syntactical
    classes that dictate how they should interact. We tart with a class of basic
    types denoted $A$, $B$, $P$, $Q$, \etc.

    If we have two types $A$ and $B$, we may construct the new type $A
    \rightarrow B$, the type of functions between $A$ and $B$. Given a variable
    $x$ of type $A$ and a term $t$ (possibly containing $x$) of type $B$, we may
    form the term $\lambda x \mathrel. t$ of type $A \rightarrow B$. If we have
    both $\function f A B$ and $t : A$ (the colon is read as ``has the type''),
    we may form $(f t) : B$.

    You might already begin to see a relation between both systems\dots\
    something we might call a \emph{correspondence}. Indeed, the Curry-Howard
    Correspondence is an observation on some kind of interdeductibility bertween
    the Simply-Typed $\lambda$-calculus and the Intuitionistic Propositional
    Logic. The rules presented above can be seen as deduction rules, and, just
    like that, we see ourselves writing
   
    \[
        \inferrule {
            \inferrule {
                \inferrule {
                    \inferrule {
                        \inferrule { } {
                            \vdash A \type
                        }
                    } {
                        x : A \vdash x : A
                    } \\
                    \inferrule {
                        \inferrule { } {
                            \vdash A \type
                        } \\
                        \inferrule { } {
                            \vdash B \type
                        }
                    } {
                        x : A \vdash B \type
                    }
                } {
                    x : A, y : B \vdash x : A
                }
            } {
                x : A \vdash \lambda y x : B \rightarrow A
            }
        } {
            \vdash \lambda x \lambda y x : A \rightarrow (B \rightarrow A)
        }
    \]

    This can be very well summarized in a theorem, paraphrasing
    \citet{SorensenLectures2006},

    \begin{theorem*}
        \begin{enumerate} It is true that

            \item If $\Gamma \vdash t : A$ is deductible in the
                $\lambda$-calculus, then $\tilde\Gamma \vdash A \true$ is
                deductible in the Propositional Logic (where $\tilde\Gamma$ is
                $\Gamma$ without variables).

            \item If $\Delta \vdash A \true$ is deductible in the Propositional
                Logic, then $\Gamma \vdash t : A$ is deductible in the
                $\lambda$-calculus for some term $t$ and $\Gamma$ having
                $\tilde\Gamma = \Delta$.

        \end{enumerate}
    \end{theorem*}

    The details of this theorem are the subject of the present work.

    \subsection{Notes on methodology}

    Some methodological particularities are to be mentioned, to the cost of some
    readers being spooked upon contact with the rest of the text. If you, the
    reader, have skimmed the rest of these pages, no doubt you will have noticed
    the presence of code listings. That's explained by two facts:

    \begin{enumerate}

        \item Our base theory is \emph{not} \textsc{ZFC}, as much as it is not
            at all classic. We're working informally in the Martin-Löf Type
            Theory, that is intuitionistic. Every argument (though there are not
            many) are \emph{constructivist} in nature, and all definitions are
            focused in pinning down \emph{structure}, not concepts.

        \item Formally, our definitions and theorems are stated in the
            type-theory which forms the basis of Idris. It is also derived from
            the Martin-Löf Type Theory. There's comprehensive material on that
            programming language online, besides the book by
            \citet{BradyDevelopment2017}. Idris is a \emph{propgramming}
            language, focused on software development, but it serves just as
            well as a formal theorem-proving enviroment.

    \end{enumerate}

    Some statements come also in an Idris version -- checked by its compiler --,
    though we may have ommited some others in the name of brevity.

    \subsection{Take a map}

    This text is structured in the following way.

    \begin{itemize}

        \item Section~\ref{sec:what-is-a-logic} talks about our formalism, and what
            is the reasoning behind its definition. It seeks to answer the
            question, as the title indicates, \emph{what is a logic?}
            Particullarly, what is a formulation of the ideia of \emph{one}
            logic that is succeptible to the Curry-Howard Correspondence.

        \item Section~\ref{sec:what-is-the-curry-howard-correspondence} uses the
            formalism developed in ~\ref{sec:what-is-a-logic} and applies it
            into formalizing the theorem.

        \item Section~\ref{sec:further} discusses possible further developments.

    \end{itemize}

    \section{What is a logic?}
    \label{sec:what-is-a-logic}

    \subsection{As for the literature}

        The natural starting point for a definitional question like the one we
        ask (``what \emph{is} the Curry-Howard Correspondence?'') is
        bibliographic revision. We hoped to find, among the miriad stablished
        definitions of what might a logic be, one that fit our bill.

        The study of logic -- which, despite the wide-spread usage of the word,
        isn't a homogeneous field -- is almost as old as western thought.
        Still, there doesn't seem to be any sort of consensus on its object of
        study. There isn't a single definition of the word ``logic''
        \citep[see][2]{SEPLogicOntology}.

        Nonetheless, there are several definitions with varying levels of
        generality. Starting with the concept of \strong{Tarski consequence
        operators} \citep{SuszkoRemarks1958}, layed out as a way to define
        logics by their notions of \emph{deductibility}, abstracted away from
        the details of the particular deductions. They stand out for being a
        ``universal'' formulation: a logic, so to speak, is determined by a set
        of data (a set of formulas and an operator) satisfying certain
        properties (like finitariety, monotonicity, reflexivity, \etc) --
        not unlike the way many structures are defined by their universal
        properties in Category Theory, only later for examples to be
        constructed.

        Another idea, more categoriacal in nature, and originated in other
        contexts, is that of \strong{institutions} -- or, for our
        proof-theoretical purposes, their cousings the $\pi$-institutions (the
        first occurence of the subject can be found in
        \citet{FiadeiroStructuring1988}). A $\pi$-institution is defined by a
        sequence of data:

        \begin{enumerate}

            \item A category $\Sig$ of signatures,

            \item A functor $\function \Sen \Sig \Set$ that, for each signature,
                constructs the set of sentences over it, and

            \item A notion $\mathcal{C}_\Sigma$ of deductibility \foreign{à lá}
                Tarski for every signature $\Sigma$ over the sentences in
                $\Sigma$, such that

            \item these data satisfy a \emph{coherence condition} over the
                signature morphisms: for $\function f \Sigma {\Sigma'}$, $\Gamma
                \subseteq \Sen\Sigma$, $\varphi \in \Sen\Sigma$

                \[
                    \Gamma \vdash_\Sigma \varphi \implies (\Sen f)[\Gamma]
                    \vdash_{\Sigma'} (\Sen f)(\varphi)
                \]

        \end{enumerate}

        It's worth noting that a \strong{$\pi$-institution morphism} is given by
        a functor between the signature categories and a natural transformation
        between the sentence functors (mediated by the morphism one), again
        satisfing some notion of coherence. Coloquially, the functor provides a
        reinterpretation of conectives, and the natural transformation a
        ``change in the form'' of the sentences.
       
        These structures can be seen as a categorical reinterpretation of Tarski
        Operators. Everything that applies to these applies point-wise to those,
        plus the additional structure of the ``translations'' given by the
        signature morphisms. A stardard reference on the subject of
        institutions, where one can get a better understanding of the workings
        of this sort of machinery, is \citet{DiaconescuInstitution2008}.

        There are even \emph{more} categorical aproaches to general definitions
        of logics. We could mention, for example,
        \citet{LambekCategorical1988}'s interpretation of higher order logics
        as categories closed under certain operations. Particularly, the idea
        that intuitionistic propositional logics correspond to cartesian-closed
        categories -- and that lambda-terms are arrows in that category.

        Finally, the least standard reference we'll considere comes from the
        famous\footnote{Maybe infamous.} n-categorical encyclopedia $n$Lab
        \citep{nLabDeductive}. A \strong{deductive system}, according to
        that definition, is ``a collection of judgements'' together with a
        ``collection of \emph{steps}'', which consist in a list of
        judgements\footnotemark -- its hypothesis -- and another judgement --
        its conclusion.

        \footnotetext{For a more comprehensive discussion on
        \strong{judgements}, see \citet{LofMeaning1996}}.
        
        \[
            \inferrule {
                J_0, \\
                \cdots \\
                J_{n-1}
            } {
                J
            }
        \]

        These ``steps'', obviously, can be composed to form deductions.
        Actually, deductions are freely generated from their composition.

        In short, we consider the following definitions and interpretations:

        \begin{description}

            \item[Tarski] The idea of logic as a ``universal'' formulation: an
                abstract structure satisfying some properties.

            \item[Institutions] Logics form a category with their translations.

            \item[Lambek] Logics are categories \emph{themselves}, or categories
                with additional structure.

            \item[Deductive systems] Deduction from proof. A logical system is
                determined by how you are allowed to construct arguments. A
                logical step is an abstract entity that allows you to jump from
                hypothesis to a conclusion.

        \end{description}

        With the goal of defining a class of systems subject to a theorem like
        Curry-Howard, we start exploring deductive systems. Let us keep in mind,
        though, that eventually we'd like to formalize the theorem between the
        already mentioned systems, with the possibility of extensions.

    \subsubsection{Deductive Systems}

        The fundamental idea of deductive systems is that deductions are
        governed by rules -- objects that authorize transitions between
        judgements. These rules, or steps, then, act over some kind of
        judgements $\mathcal{J} : \Type$. There are pairs, so to speak, that
        indicate when a conclusion $J : \mathcal{J}$ can be made from hypothesis
        $J_1, \cdots, J_{n-1} : \mathcal{J}$. Formally,

        \begin{definition}[Step/Inference]

            An \strong{inference} over a base type $\mathcal{J} : \Type$ is a
            triple containing the indices of its hypothesis\footnotemark, the
            hypothesis themselves, and a conclusion.

            \begin{code*}{firstnumber=14,label=idr/DeductiveSystems/DeductiveSystems.idr}
                record Inference j where
                  constructor MkInference
                  labels : Type
                  hypothesis : labels -> j
                  conclusion : j
            \end{code*}

        \end{definition}

        \footnotetext{Note our choice to use indexed families instead of a list,
        or a set. That offers, in general, quality of life for the theory's
        developer. Your are welcomed to try other formalisms in any case of
        disagreement.}

        A deductive system is determined not only from its judgements or logical
        steps, but also by \emph{which} steps ever are authorized. We call that
        extra structure a \strong{ruleset}.

        \begin{definition}[Ruleset]

            A \strong{ruleset} is an indexed family of inferences.

            \begin{code*}{firstnumber=26,label=idr/DeductiveSystems/DeductiveSystems.idr}
                Ruleset : Type -> Type -> Type
                Ruleset l j = (l -> Inference j)
            \end{code*}

        \end{definition}

        Now, we define what we actually want to achieve with deductive systems.
        In a way, we have defined up to here the \emph{signature} of a deductive
        system, but not its deductions. So for every judgement, we attribute it
        the \strong{type of its proofs}.

        \begin{code*}{firstnumber=36,label=idr/DeductiveSystems/DeductiveSystems.idr}
            data Deduction : Ruleset l j -> j -> Type where
        \end{code*}

        Constructing a proof, given a type for the judgements and a
        ruleset over it, involves given a step (that's \emph{in} the ruleset)
        and proofs of its hypothesis. That way, we obtain a proof of its
        conclusion.

        That is, proofs are build inductively, or freely, from the application
        of the following function.

        \begin{code*}{firstnumber=46,label=idr/DeductiveSystems/DeductiveSystems.idr}
              Infer : {j : Type} ->                  -- Given a base type J,
                      {r : Ruleset l j} ->           -- a ruleset R,
                      {i : l} ->                     -- and an inference r,
                      (h : labels (r i)) -> Deduction r (hypothesis (r i) h) ->
                      Deduction r (conclusion (r i)) -- we get a proof of the conclusion
        \end{code*}

    \subsection{Relations}

        We may test definitions by comparing them to similar, better stablished
        ones. As we work with deductive systems, inspired in an entry to a
        non-standard encyclopedia, we might ask ourselfs: how does this compare
        to other formalisms we've mentioned? Let's start with Tarski operators
        -- or their equivalent counterparts, deductive relations.

        Right off the bat we observe that deductive systems, thus abstractly
        defined, is not a theory of \emph{inference}, just as much as it is a
        theory of \emph{proofs} -- a proof-theory. After all, we've defined a
        judgement's type of proofs, with no space for hypothesis. To this
        purpose, we may \emph{refine} the type of judgements. Instead of a
        generic type, we consider a type of \strong{hypothetical} judgements
        over a base type.

        \begin{definition}[Inferência]

            A \strong{inference} over a type of judgements $\Judgement$ is an
            object of the form

            \[
                J_0, \cdots, J_{n-1} \vdash J
            \]

            where $J_i, J : \Judgement$. We denote the type of these objects as
            $\Inference\;\Judgement$.

        \end{definition}

        \begin{definition}

            A \strong{hypothetical ruleset} over a type $\Judgement$ is a
            ruleset over $\Inference\;\Judgement$.

        \end{definition}

        It's now evident the translation we seek: a hypothetical deductive
        system determines a deductive relation, where we say that a certain
        inference is valid under that relation if there's proof of it. That is,
        we say $\Gamma \succ J$ if there's proof of $\Gamma \vdash
        J$\footnotemark.

        \footnotetext{We denote deductive relations by $\succ$ and hypothetical
        judgements by $\vdash$, to avoid misunderstandings.}

        But alas, here we find a problem. On the one hand, the definition of an
        inference construes its hypothesis as a finite list of judgements, which
        is usual and sometimes \emph{necessary}\footnotemark in various type
        theories. On the other hand, deductive relations are most commonly
        instituted between \emph{subsets} of the type of judgements and
        judgements themselves.

        \footnotetext{In dependent type theories, contexts (hypothesis) must be
        \emph{ordered} because later types might depend on variables declared to
        be of earlier ones. Finiteness is a common property, as we often
        construct proofs as necessarily finite objects.}

        So how might we correct such disparity? We could hope to define
        translations from one kind of judgement to the other. As the base type
        is the same, we've only been left with \emph{transforming}: literally
        changing the \emph{form} of one kind of judgement into the other. Here
        we might notice
       
        \begin{mathpar}
            \Gamma \vdash J \and \Delta \succ J
        \end{mathpar}

        that the relevant difference is just in the form of the hypothesis:
        finite lists on the left, and subsets on the right. That's the
        transformation we want to achieve.

        The seemingly most natural way to do it is to take the \strong{range} of
        a list to form a set, and to take finite lists from a set the other way
        around. That is, denoting by ``$\Gamma \vdash J$'' the phrase ``there's
        proof $\rho : \Deduction_R \Gamma \vdash J$'',

        \begin{definition}

            The deductive relation $\succ_\vdash$ determined by a deductive
            system $\vdash$ is given by

            \begin{center}
                $\Delta \succ J$ if there is $\Delta$-sequence $\Gamma$ such
                that $\Gamma \vdash J$
            \end{center}

            Where a $\Delta$-sequence is a finite list with its range a subset
            of $\Delta$.

            The deductive system $\vdash_\succ$ determined by a deductive
            relation $\succ$ is given by

            \begin{center}
                $\Gamma \vdash J$ if $\img(\Gamma) \succ J$
            \end{center}

        \end{definition}

        A keen reader might have noticed that these definitions seem to define
        some kind of adjection between the classes of deductive systems and
        relations. Indeed,

        \begin{proposition}

            It is true that

            \begin{align*}
                \Gamma \vdash_{\succ_\vdash} J &\iff \Gamma \vdash J \\
                \Delta \succ_{\vdash_\succ} J &\implies \Delta \succ J \quad
                \text{if $\succ$ is finitary}\\
                \Delta \succ_{\vdash_\succ} J &\impliedby \Delta \succ J \quad
                \text{if $\succ$ is monotonic}
            \end{align*}

            A relation is finitary if from an inference we may find another that
            stands between a finite subset of the original hypothesis and the
            same conclusion. A system is monotonic if we may add hypothesis
            whilst preserving the inference.
                
            \begin{proof}

                Starting with the first one,

                \begin{align*}
                    \Gamma \vdash_{\succ_\vdash} J 
                    &\iff \img\Gamma \succ_\vdash J \\
                    &\iff \text{exists $\img\Gamma$-sequence $\Gamma'$ such that
                    $\Gamma' \vdash J$}
                \end{align*}

                But $\Gamma$ is a $\img\Gamma$-sequence and $\Gamma \vdash J$.

                For the later two implications,

                \begin{align*}
                    \Delta \succ_{\vdash_\succ} J
                    &\iff \text{exists $\Delta$-sequence $\Gamma$ such that $\Gamma
                    \vdash_\succ J$} \\
                    &\iff \text{exists $\Delta$-sequence $\Gamma$ such that
                    $\img\Gamma \succ J$}
                \end{align*}

                But the image of a $\Delta$-sequence is exactly a finite subset
                of $\Delta$. So we conclude that

                \begin{itemize}

                    \item If $\succ$ is finitary and $\Delta \succ J$, we take
                        finite $\Delta' \subseteq \Delta$ such that $\Delta'
                        \succ J$. Any enumeration $\Gamma$ of $\Delta'$ is such
                        that

                        \begin{align*}
                            \Delta = \img\Gamma &\succ J \\
                            \Gamma &\vdash_\succ J
                        \end{align*}

                        But $\Gamma$ is a $\Delta$-sequence.

                    \item If $\succ$ is monotonic, and $\Delta
                        \succ_{\vdash_\succ} J$, then there exists a
                        $\Delta$-sequence $\Gamma$ such that $\Gamma
                        \vdash_\succ J$ -- \ie, $\img\Gamma \succ J$. But
                        $\img\Gamma \subseteq \Delta$. So $\Delta \succ J$.

                \end{itemize}

            \end{proof}

        \end{proposition}

        We've stumbled upon a concept that's reminicent of a \emph{morphism}
        between these two logical systems. But instead of exchanging signatures,
        as is normally the case, we've \emph{transformed} them -- literally,
        changed their form.

        That's familiar enough: it's analogous to a $\pi$-institution morphism.
        Though we've kept fixed a signature, the change-in-form behaves as some
        kind of natural transformation between the ``sentence functors'' of the
        deductive systems and relations.

        A troubling inconsistency of this interpretation, however, is the fact
        that deductive systems aren't as ``abstract'' as deductive relations.
        We've defined them explicitly, demanding data that informs us
        \emph{exactly} how to construct proofs. They talk about
        \emph{proof}, not \emph{provability}.

        This can, however, be fixed.

    \subsubsection{p-institutions}

        We start our definitive formalization with the notion of a
        \strong{predicate}, the ``abstract'' version of a deductive system.

        \begin{definition}[Predicate]

            A \strong{predicate} over a type is just an indexed family over it.
            Under the propositions-as-types interpretation -- which, besides
            being formalized by us provides us with a foundational principle in
            our methodology --, they may be understood as a function that
            associates to every element of a base type a proposition -- \ergo, a
            predicate.

            In another, complimentary reading, a predicate pairs every
            inhabitant of its universe of discourse with the type of its
            \emph{proofs}. So that's the \emph{proofful} interpretation.

            \begin{code*}{firstnumber=19,label=idr/Predicates/Predicates.idr}
                Predicate : Type -> PointedFunctor TypeCat -> Type
                Predicate a f = (mapObj (functor f) a) -> Type
            \end{code*}

        \end{definition}

        We define also an appropriate notion of morphism between these objects:

        \begin{definition}[Predicate Morphism]

            A \strong{predicate morphism} is a function between base types that
            preserves provability.

            In more detail, a morphism is a function that lifts to the type of
            proofs. \ie, given predicates $\function \pi a \universe$ and
            $\function \tau b \universe$, a morphism is a function $\function f
            a b$ such that, for every $j : a$, there's a function $\function
            {\varepsilon_j} {\pi(j)} {\tau(f(j))}$.

            \begin{code*}{firstnumber=33,label=idr/Predicates/Predicates.idr}
                record PredicateMorphism (a : Type)
                                         (b : Type)
                                         (f : PointedFunctor TypeCat)
                                         (g : PointedFunctor TypeCat)
                                         (p : Predicate a f)
                                         (q : Predicate b g) where
                  constructor MkPredicateMorphism
                  Translate : a -> b
                  Transform : PointedNaturalTransformation TypeCat f g
                  Transport : {x : mapObj (functor f) a} -> p x -> q (apply Transform Translate x)
            \end{code*}

        \end{definition}

        \begin{definition}[p-Institutions]

            A \strong{p-institution} is composed of a signature category $\Sig$,
            a functor $\function \Sen \Sig \Type$ and a family of functions

            \[
                \pi : \prod_{a : \Sig} \mathrm{Predicate}\paren{
                    \Sen\,a
                }
            \]

            that respect a coherence condition

            \[
                \prod_{a,b : \Sig}
                \prod_{\function f a b}
                \prod_{j : \Sen\,a}
                \pi_a\,j \longrightarrow
                \pi_b\,(\Sen\,f)\,j
            \]

        \end{definition}

        \begin{definition}[p-Institution Morphism]

            A \strong{p-institution morhpism} between $(\Sig, \Sen, \pi)$ and
            $(\Sig', \Sen', \pi')$ is determined by a functor $\function F \Sig
            {\Sig'}$ and a natural transformation $\function \mu \Sen {\Sen'
            \circ F}$ such that

            \[
                \prod_{a : \Sig}
                \prod_{j : \Sen\,a}
                \pi\,j \longrightarrow
                \pi' (\mu\,j)
            \]

        \end{definition}

    \subsubsection{Deductive systems as p-institutions}

        We've defined p-institutions in order to have an abstract notion of
        deductive systems. A particular deductive system correspondes to a
        predicate, but a p-institution is a structured collection of predicates.
        So, in order to see deductive systems as p-institutions, we have to
        consider structure collections of them. Particularly, a kind of indexed
        family.

        We start with a category $\Sig$ and a functor $\function \Sen \Sig
        \Type$. But now, instead of attributing a predicate to every signature,
        we attribute a ruleset. For every $\Sigma : \Sig$ we attribute an
        indexing family $l_\Sigma : \Type$ and a ruleset $R_\Sigma : \Ruleset \,
        l_\Sigma \, (\Sen \Sigma)$. Besides, for every arrow $\function f \Sigma
        {\Sigma'}$ a ruleset morphism $\function {R_f} {R_\Sigma} {R_\Sigma'}$.

        Given the data $(\Sig, \Sen, l, R)$ we may define the p-institution
        given by

        \[
            \paren{
                \Sig,
                \Sen,
                \pi
            }
        \]

        where

        \[
            \pi \, a \, J = \Deduction \, R_a \, J
        \]

        as expected\footnotemark.

        \footnotetext{It's interesting to consider that a pair $(\Sig, \Sen)$
        might admit several notions of provability, made concrete in the form of
        several rulesets over it. Considering as arrows between these pairs
        functors and natural transformations as previously defined, could wi set
        up a $\pi$ or p-institution of the ``deductive structures over $(\Sig,
        \Sen)$? We leave you with that.}

\section{What is the Curry-Howard Correspondence?}
\label{sec:what-is-the-curry-howard-correspondence}

    \subsection{Propositional Logic}

        The tradicional Curry-Howard Theorem, as presented in
        \citet{SorensenLectures2006}, is a specfic case on a distinguished
        result about a particular morphism between two p-institutions. One of
        them, which we explore first, it that of the so called
        \strong{propositional logics}.

        We start, as usual, with signatures.

        \begin{definition}[Propositional Signature]

            A \strong{propositional signature} is determined by an indexed
            family of arities, rules and a type for variables. That is, a tuple
            $(\connective, \arity, \rules, \headers, \var)$ where

            \begin{itemize}

                \item $\connective : \Type$, is a type whose inhabitants are
                    called \strong{connectives}.

                \item $\function \arity \connective \Naturals$, an arity
                    function.

                \item For connective $c : \connective$, a family $\headers_c =
                    \unitary{h_r}_{r : \rules(c)}$ (where $function \rules
                    \connective \Type$) of \strong{headers}: lists containing
                    the elements $\true$ or $\prop$, whose lenght is the
                    connectives arity $\arity c$.

                \item A type $\var : \Type$ for variables.

            \end{itemize}
            
            \begin{code*}{firstnumber=17,label=idr/DeductiveSystems/PropositionalDeductiveSystems.idr}
                record Signature where
                  constructor MkSignature
                  connective : Type
                  arity : connective -> Nat
                  rule : connective -> Type
                  header : (c : connective) -> (rule c) -> Header (arity c)
                  var : Type
            \end{code*}

        \end{definition}

        A signature, in the context of anything-institutions, can be understood
        as the \emph{most concise} data that distinguishes the different logics
        inside the same institution-like structure. In the case of propositional
        logics, that means their connectives, arities (which, by themselves,
        already determine all propositions) besides the rules governing them
        (which determine the notion of proof).

        \begin{definition}[Propositional Signature Morphism]

            A \strong{propositional signature morphism} is determined by a
            function between the connectives that preserves arity and that's
            lifted to a function between rules, besides a separate function
            between variables.

            That is, given signatures $(\connective, \arity, \rules, \headers,
            \var)$ and $(\connective', \arity', \rules', \headers', \var')$, a
            function $\function f \connective {\connective'}$ such that $\forall
            c \in \connective \mathrel. \arity c = \arity'f(c)$, and functions
            $\function {g_c} {\rules(c)} {\rules'(f(c))}$ and $\function {h}
            {\var} {\var'}$.

        \end{definition}

        \begin{remark}

            The given definition for morphisms it, for the most part, too
            strict. It is, however, simpler than the alternatives -- and, for
            our purposes, it doesn't make any difference. We adopt it in the
            name of simplicity.

        \end{remark}

        \begin{definition}[Propositions over a signature]

            The \strong{propositions} over a signature are build inductively:

            \begin{itemize}

                \item A variable is a proposition, and

                \item A connective with arity $n$ applied to $n$ propositions is
                    a proposition.

            \end{itemize}

            \begin{code*}{firstnumber=28,label=idr/DeductiveSystems/PropositionalDeductiveSystems.idr}
                data Proposition : Signature -> Type where
                  Atomic : var s -> Proposition s
                  Apply : {s : Signature} ->
                          (c : connective s) ->
                          Vect (arity s c) (Proposition s) ->
                          Proposition s
            \end{code*}

        \end{definition}

        But our logic's judgements are not propositions: they are \emph{about}
        propositions, under certain hypothesis \citep[see][]{LofMeaning1996}.

        \begin{definition}[Propositional Judgements]

            The \strong{propositional judgements} are of the form

            \[
                P_0, \cdots, P_n \vdash P\,\mathit{adjective}
            \]

            where $P_i$ and $P$ are propositions and $\mathit{adjective}$ is
            either $\true$ or $\prop$.

            \begin{code*}{firstnumber=52,label=idr/DeductiveSystems/PropositionalDeductiveSystems.idr}
                Judgement : Signature -> Type
                Judgement s = (Context s, Lof (Proposition s))
            \end{code*}

            The object to the right of $\vdash$ is a \strong{atomic judgement}.

            \begin{code*}{firstnumber=38,label=idr/DeductiveSystems/PropositionalDeductiveSystems.idr}
                Lof : Type -> Type
                Lof j = (j, LofAdjective)
            \end{code*}

            The object to its left is a \strong{context} -- a finite list of
            propositions.

            \begin{code*}{firstnumber=46,label=idr/DeductiveSystems/PropositionalDeductiveSystems.idr}
                Context : Signature -> Type
                Context s = List (Proposition s)
            \end{code*}

        \end{definition}

        Finally, we come to the rules.

        The rules of propositional logic split into a handful of species: axioms
        and structural, syntactical, introduction and elimination rules. Given a
        signature $\Sigma = (\connective, \arity, \rules, \headers, \var)$,

        \begin{description}

            \item[Axioms] For every propositional variable $p : \var$, postulate
                that

                \[
                    \inferrule { } {
                        \vdash p \prop
                    } \quad \textsc{Var}_p
                \]

            \item[Structural] For every context $\Gamma$, proposition $P$ and
                judgement $J$, allow it to be that

                \begin{mathpar}
                    \inferrule {
                        \Gamma \vdash P \prop
                    } {
                        \Gamma, P \vdash P \true
                    } \quad \textsc{Self} \and
                    \inferrule {
                        \Gamma \vdash P \prop \\
                        \Gamma \vdash J
                    } {
                        \Gamma, P \vdash J
                    } \quad \textsc{Weak}
                \end{mathpar}

            \item[Syntactical] Every connective $c : \connective$ has a
                \strong{formation rule}: given propositions $P_0, \cdots,
                P_{\arity c}$ and context $\Gamma$,

                \[
                    \inferrule {
                        \Gamma \vdash P_0 \prop \\
                        \cdots \\
                        \Gamma \vdash P_{\arity c} \prop
                    } {
                        \Gamma \vdash c(P_0, \cdots, P_{\arity c}) \prop
                    } \quad \textsc{Form}_c
                \]

            \item[Introduction] For every connective $c : \connective$ and rule
                $r : \rules(c)$, for propositions $P_0, \cdots, P_{\arity c}$
                and context $\Gamma$,

                \[
                    \inferrule {
                        \Gamma \vdash P_0 \headers_c(r)[0] \\
                        \cdots \\
                        \Gamma \vdash P_{\arity c} \headers_c(r)[\arity c]
                    } {
                        \Gamma \vdash c(P_0, \cdots, P_{\arity c}) \true
                    } \quad \textsc{Intro}_c^r
                \]

            \item[Elimination] For every connective $c$, an \strong{elimination
                rule}

                \[
                    \inferrule {
                        \Gamma \vdash c(P_0, \cdots, P_{\arity c}) \\
                        \Gamma, P_{r_11}, \cdots, P_{r_1k_1} \vdash C \true \\
                        \cdots \\
                        \Gamma, P_{r_l1}, \cdots, P_{r_lk_l} \vdash C \true
                    } {
                        \Gamma \vdash Q
                    } \quad \textsc{Elim}_c
                \]

                for context $\Gamma$, propositions $P_0, \cdots, P_{\arity c}$
                and proposition $C$. The propositions $P_{ri}$ are chosen among
                the $P_j$, where $r$ varies over the rules and $i$ over $r$'s
                hypothesis that have the adjective $\true$.

        \end{description}

        Let's look at an example: suppose a signature has among its connectives
        a  $\lor : \connective$, whose arity is $\arity \lor = 2$, whose rules
        are $\Left : \rules_\lor$ e $\Right : \rules_\lor$, and headers

        \begin{align*}
            \headers_c(\Left) &= [{\true}, {\prop}] \\
            \headers_c(\Right) &= [{\prop}, {\true}]
        \end{align*}

        So we have the following rules

        \begin{mathpar}
            \inferrule {
                \Gamma \vdash A \prop \\
                \Gamma \vdash B \prop
            } {
                \Gamma \vdash A \lor B \prop
            } \quad \textsc{Form}_\lor         \\
            \inferrule {
                \Gamma \vdash A \true \\
                \Gamma \vdash B \prop
            } {
                \Gamma \vdash A \lor B \true
            } \quad \textsc{Intro}_\lor^\Left  \and
            \inferrule {
                \Gamma \vdash A \prop \\
                \Gamma \vdash B \true
            } {
                \Gamma \vdash A \lor B \true
            } \quad \textsc{Intro}_\lor^\Right \\
            \inferrule {
                \Gamma \vdash A \lor B \true \\
                \Gamma, A \vdash C \true \\
                \Gamma, B \vdash C \true
            } {
                \Gamma \vdash C \true
            } \quad \textsc{Elim}_\lor

        \end{mathpar}

        Note that despite not expressing it explicitly, we have specified the
        index besides the rules themselves. The index of a rule consists on the
        name $\textsc{Form}_c$, $\textsc{Intro}_c^r$, \etc, besides the data
        introduced as in ``for every context $\Gamma$ and propositions
        \omissis''. The rules themselves are given by their diagramatical
        representations.

        That said, we can take a look at a sample demonstration. We leave it as
        an exercise to fill-in the names of the rules, complete the missing
        parts or even \emph{redoing} the proof-tree without looking at it.

        \[ \footnotesize
            \inferrule {
                \inferrule {
                    \inferrule {
                        \inferrule {
                            \inferrule { } {
                                \vdash B \prop
                            } \\
                            \inferrule { } {
                                \vdash A \prop
                            }
                        } {
                            \vdash A \lor B \prop
                        }
                    } {
                        A \lor B \vdash A \lor B \true
                    } \\
                    \inferrule {
                        \inferrule {
                            \inferrule {
                                \inferrule { }{
                                    \vdash A \prop
                                } \\
                                \inferrule { }{
                                    \vdash B \prop
                                }
                            } {
                                \vdash A \lor B \prop
                            } \\
                            \inferrule { } {
                                \vdash A \prop
                            }
                        } {
                            A \lor B \vdash A \prop
                        }
                    } {
                        A \lor B, A \vdash A \true
                    } \\
                    \inferrule {
                        \inferrule {
                            \omissis
                        } {
                            A \lor B \vdash B \prop
                        } \\
                        \inferrule {
                            \omissis
                        } {
                            A \lor B \vdash A \true
                        }
                    } {
                        A \lor B, B \vdash A \true
                    }
                } {
                    A \lor B \vdash A \true
                }
            } {
                A \lor B \vdash B \lor A \true
            }
        \]

        Note, finally, how the proof-tree has certain ``transition'' nodes
        between judgements of type $\true$ and of type $\prop$. When does that
        transition happen?

        \thematicbreak

        Formally, the rules' indices form an inductive type, as follows (just
        the introduction and elimination cases).

        \begin{code*}{firstnumber=61,label=idr/DeductiveSystems/PropositionalDeductiveSystems.idr}
            data RuleName : Signature -> Type where
              Intro : (c : connective s) ->
                      (r : rule s c) ->
                      (ctx : Context s) ->
                      Vect (arity s c) (Proposition s) ->
                      RuleName s
              Elim : (c : connective s) ->
                     (ctx : Context s) ->
                     Vect (arity s c) (Proposition s) ->
                     Proposition s ->
                     RuleName s
        \end{code*}

        The \emph{ruleset} itself is an indexed famoly \idris{RuleName}
        $\Sigma$, a function pairing

        \begin{itemize}

            \item The name \idris{Intro c r ctx ps} to the rule

                \[
                    \inferrule {
                        \idris{zip ps (header s c r)}
                    } {
                        \idris{Apply c ps true}
                    }
                \]

            \item And the name \idris{Elim c ctx ps p} to

                \[
                    \inferrule {
                        \overbrace {
                            \idris{ctx} \vdash \idris{Apply c ps true}
                        }^{\idris{Left ()}} \\
                        \overbrace {
                            \idris{ctx ++ filterWith ps (rule s c r) IsIsTrue} \vdash \idris{p}
                        }^{\idris{Right r}}
                    } {
                        \idris{ctx} \vdash \idris{p}
                    }
                \]

                The function \idris{filterWith} takes from the vector \idris{ps}
                those elements which are paired with elements of the header
                \idris{rule s c r} that satisfy the boolean predicate
                \idris{IsIsTrue} -- those propositions that, for that header,
                need to be true in the application of that rule --, e puts them
                together in a list that's appended to the context.

        \end{itemize}

        The previous example -- the disjunction -- would have an elimination
        rule as in

        \[
            \inferrule {
                \idris{ctx} \vdash \idris{(Apply Or [p1, p2]) true} \\
                \idris{ctx ++ [p1]} \vdash \idris{p true} \\
                \idris{ctx ++ [p2]} \vdash \idris{p true}
            } {
                \idris{ctx} \vdash \idris{p true}
            }
        \]

    \subsection{$\lambda$-calculus}

        \begin{description}

            \item[Signatures] The signatures of the $\lambda$-calculus are
                identical to the propositional logic's ones, except for the
                addition of a single type for individual variables. That is,

                \[
                    (\connective, \arity, \rules, \headers, \var_\ty, \var_\tm)
                \]

                The morphisms are extended in the obvious way.

            \item[Judgements] The judgements of the $\lambda$-calculus over
                $\Sigma = (\connective, \arity, \rules, \headers, \var)$ are of
                the form

                \[
                    v_0 : P_0, \cdots, v_{n-1} : P_{n-1} \vdash \mathcal{J}
                \]

                where $v_i : \var_\tm$, $P_i$ are propositions over $\Sigma$ and
                $J$ has one of the forms
                
                \begin{mathpar}
                    P \prop \and
                    t : P
                \end{mathpar}

                where $P$ is a proposition and $t$ is a \strong{term} over
                $\Sigma$:

                \begin{itemize}

                    \item Every individual variable is a term.

                    \item For every connective $c : \connective$, and every rule
                        $r : \rules$, there's a term \strong{constructor}
                        $\lambda_c^r$ that takes as many terms as there are
                        hypothesis of the form $t : A$; and returns a term.

                    \item For every connective $c : \connective$, there's an
                        \strong{eliminator} $\application_c$ that takes, for
                        every rule $r : \rules$, as many variables as there are
                        hypothesis of the form $t : \varphi$; and a term, and
                        teruns a term. \eg,

                        \begin{mathpar}
                            \inferrule {
                                \Gamma \vdash a : A \\
                                \Gamma \vdash B \type
                            } {
                                \Gamma \vdash \lambda_\lor^L(a) : A \lor B
                            } \and
                            \inferrule {
                                \Gamma \vdash A \type \\
                                \Gamma \vdash b : B
                            } {
                                \Gamma \vdash \lambda_\lor^R(b) : A \lor B
                            } \and
                            \longmapsto \and
                            \inferrule {
                                \Gamma \vdash p : A \lor B \\
                                \Gamma, x : A \vdash c_a : C \\
                                \Gamma, y : B \vdash c_b : C
                            } {
                                \Gamma \vdash \application_\lor(p, x, c_a, y, c_b) : C
                            }
                        \end{mathpar}

                        In the example, the constructors are usually denoted

                        \begin{align*}
                            \lambda_\lor^L(a) &\longmapsto \InL(a) \\
                            \lambda_\lor^R(b) &\longmapsto \InR(b) \\
                            \varepsilon_c(p, x, c_a, y, c_b) &\longmapsto
                            \caseof{p}{x}{c_a}{y}{c_b}
                        \end{align*}

                        Note that the instances of the symbols $\InL$ and $\InR$
                        in the eliminator are simply \emph{asthetic} -- not
                        really applications of the constructors.

                \end{itemize}

            \item[Rules]

                \begin{description}

                    \item[Axiom] For every propositional variable $p : \var_\ty$
                        and context $\Gamma$, an axiom

                        \[
                            \inferrule {
                            } {
                                \Gamma \vdash p \type
                            } \quad \textsc{Var}(\Gamma,p)
                        \]

                    \item[Structural] For every context $\Gamma$, variable $x$
                        that does not occur in $\Gamma$, proposition $P$ and
                        judgement $J$,
                        
                        \begin{mathpar}
                            \inferrule {
                                \Gamma \vdash P \type
                            } {
                                \Gamma, x : P \vdash x : P
                            } \quad \textsc{Self}^\lambda \and
                            \inferrule {
                                \Gamma \vdash P \type \\
                                \Gamma \vdash J
                            } {
                                \Gamma, x : P \vdash J
                            } \quad \textsc{Weak}^\lambda
                        \end{mathpar}

                    \item[Syntactical] For every connective $c : \connective$,
                        context $\Gamma$, proposition $P_i$, $0 \leq i < \arity
                        c$,

                        \[
                            \inferrule {
                                \Gamma \vdash P_0 \type \\
                                \cdots \\
                                \Gamma \vdash P_{n-1} \type
                            } {
                                \Gamma \vdash c(P_0, \cdots, P_{n-1}) \type
                            } \quad \textsc{Form}
                        \]

                    \item[Introduction] For every connective $c : \connective$
                        and rule $r : \rules(c)$, for propositions $P_0, \cdots,
                        P_{\arity c}$ and as many terms $t_i$ as there are
                        entries in $r$'s header with the value $\true$,

                        \[
                            \inferrule {
                                \Gamma \vdash J_0 \\
                                \cdots
                                \Gamma \vdash J_{\arity c}
                            } {
                                \Gamma \vdash \lambda_c^r(\mathbf{t}) : c(P_0,
                                \cdots, P_{\arity c})
                            }
                        \]

                        where $J_i$ is $P_i \type$ if $\headers_c(r) = \type$
                        and $t_i : P_i$ if $\headers_c(r) = \true$.

                    \item[Elimination] For every connective $c$, an elimination
                        rule

                        \[
                            \inferrule {
                                \Gamma \vdash p : c(P_0, \cdots, P_{\arity c}) \\
                                \Gamma, x_{r_11} : P_{r_11}, \cdots, x_{r_1k_1} : P_{r_1k_1} \vdash t_1 : C \\
                                \cdots \\
                                \Gamma, x_{r_l1} : P_{r_l1}, \cdots, x_{r_1k_1} : P_{r_lk_l} \vdash t_l : C \\
                            } {
                                \Gamma \vdash \varepsilon_c(x_{r_11}, \cdots, x_{r_1k_1}, t_1, \cdots, t_l)
                            }
                        \]

                        still standing the same remarks as for the Propositional
                        Logic eliminators\footnotemark.

                        \footnotetext{Yes, the syntax seems dense -- because it
                        is. These rules work much better in syntaxes as that of
                        Idris, in which the use of arbitrary functions
                        substitues the conventions of mathematical writing.}

                \end{description}
                
        \end{description}

\subsection{The traditional correspondece, revisited}

    Having define both systems -- the propositional logic $(\Sig^\Prop,
    \Sen^\Prop, \pi^\Prop)$ and the $\lambda$-calculus $(\Sig^\lambda,
    \Sen^\lambda, \pi^\lambda)$ --, we strive to relate them. There's an obvious
    way: a ``forgetful'' morphism. The judgements in the $\lambda$-calculus
    carry more information, in some sense, than those of Propositional Logic. We
    may, then, discard that information.

    We define a morphism $\mathcal{T}_{\lambda\to\Prop} = (T, \alpha, \beta)$
    between the p-institutions defined from both parametrized deductive systems.

    \begin{itemize}

        \item The functor between signature-categories is a forgetful one.

            \[
                T : (\connective, \arity, \rules, \headers, \var_\ty, \var_\tm)
                \longmapsto
                (\connective, \arity, \rules, \headers, \var_\ty)
            \]

        \item The natural transformation, for its part, is given by the
            following rules: for $\Sigma : \Sig$, we define

            \[
                \function {\alpha_\Sigma} {\Sen^\Prop\Sigma} {\Sen^\lambda T(\Sigma)}
            \]

            \begin{align*}
                x_0 : P_0, \cdots, x_{n-1} : P_{n-1} &\vdash t : P & &\longmapsto &
                P_0, \cdots, P_{n-1} &\vdash P \true \\
                x_0 : P_0, \cdots, x_{n-1} : P_{n-1} &\vdash P \type & &\longmapsto &
                P_0, \cdots, P_{n-1} &\vdash P \prop
            \end{align*}

        \item The proof-translation $\beta$ is given rule-by-ryle: for $\Sigma :
            \Sig$ and $\varphi : \Sen^\Prop\Sigma$,

            \[
                \function {\beta_{\Sigma, \varphi}} {\pi^\Prop_\Sigma\varphi} {\pi^\lambda_{T(\Sigma)}\beta_\Sigma\varphi}
            \]

            Writing $\equivclass{\cdot}$ for $\beta_{\Sigma, \varphi}(\cdot)$,

            \begin{align*}
                \equivclass{\textsc  {Var}^\lambda_p}
                &= \textsc{Var}_p \\
                \equivclass{\textsc {Self}^\lambda(\Ded)}
                &= \textsc{Self}(\equivclass{\Ded}) \\
                \equivclass{\textsc {Weak}^\lambda(\Ded, \Ded')}
                &= \textsc{Weak}(\equivclass{\Ded}, \equivclass{\Ded'} \\
                \equivclass{\textsc {Form}^\lambda_c(\Ded_0, \cdots, \Ded_{\arity c})}
                &= \textsc{Form}(\equivclass{\Ded_0}, \cdots, \equivclass{\Ded_{\arity c}}) \\
                \equivclass{\textsc{Intro}^{\lambda,r}_c(\Ded_0, \cdots, \Ded_{\arity c})}
                &= \textsc{Intro}(\equivclass{\Ded_0}, \cdots, \equivclass{\Ded_{\arity c}}) \\
                \equivclass{\textsc {Elim}^\lambda_c(\Ded, \Ded_{r_0}, \cdots, \Ded_{r_l})}
                &= \textsc{Elim}(\equivclass{\Ded}, \equivclass{\Ded_{r_0}}, \cdots, \equivclass{\Ded_{r_l}})
            \end{align*}

    \end{itemize}

    Inside our formalism, with all such preliminaries, the Curry-Howard
    Correspondence takes its form as a short and elegant statement.
    
    \begin{theorem}[Powered-Up Curry-Howard] $\mathcal{T}_{\lambda\to\Prop}$ is
        surjective in the space of all proofs.

        In details, given, in the propositional logic,
        
        \begin{itemize}

            \item a signature $\Sigma_\Prop : \Sig^\Prop$,

            \item a judgement $J_\Prop : \Sen^\Prop\Sigma_\Prop$, and
                
            \item a proof $\rho_\Prop : \pi^\Prop_\Sigma J$,
        
        \end{itemize}

        there exists, in the $\lambda$-calculus,
        
        \begin{itemize}

            \item a signature $\Sigma_\lambda : \Sig_\lambda$,

            \item a judgement $J_\lambda : \Sen^\lambda_\Sigma\Sigma_\lambda$,
                and
                
            \item a proof $\rho_\lambda : \pi^\lambda_\Sigma J$,

        \end{itemize}

        such that

        \begin{itemize}

            \item $T(\Sigma_\lambda) = \Sigma_\Prop$,

            \item $\alpha_{\Sigma_\lambda}(J_\lambda) = J_\Prop$ e

            \item $\beta_{\Sigma_\lambda J_\lambda} (\rho_\lambda) = \rho_\Prop$

        \end{itemize}

        \begin{proof}

            Take $(\Sigma, J, \rho$ in the space-of-proofs in the prpositional
            logic, where $\Sigma = (\connective, \arity, \rules, \headers, \var,
            \Naturals)$ and $J = P_0, \cdots, P_{n-1} \vdash P \true/{\prop}$.
            There are four steps for our demonstration: to build a signature, to
            build a judgement over that a signature, to prove it and show that
            all that data hits our target under $\mathcal{T}_{\lambda\to\prop}$.

            The signature is rather trivial: just add to $\Sigma$ a
            \emph{sufficient} type for variables, like $\Naturals$ -- \ie,
            $(\connective, \arity, \rules, \headers, \var, \Naturals)$.

            The new judgement is a bit trickier, but it's nothing more than an
            extension of the proof given by \citet{SorensenLectures2006}. We
            construct, from the proof $\rho$, a term $t$ over distinct variables
            $x_0, \cdots, x_{n-1}$\footnotemark.

            \footnotetext{In the demonstration, ``it is convenient to choose
            $\Gamma = \set{(x_\varphi : \varphi)}{\varphi \in \Delta}$, where
            $x_\varphi$ are distinct variables''
            \citep[p.~78]{SorensenLectures2006}.}

            \[
                x_0 : P_0, \cdots, x_{n-1} : P_{n-1} \vdash t : P
            \]

            or simply to consider

            \[
                x_0 : P_0, \cdots, x_{n-1} : P_{n-1} \vdash P \type
            \]

            depending on the adjective.

            The term $t$ is built by induction in $\rho$, together with the
            proof of its typing. We omit the later for simplicity, as it's just
            a matter of adding a superscript $\lambda$ to the proofs. Anyway, it
            must be the case that $\rho$ is of one of the forms

            \begin{description}

                \item[$\textsc{Var}_p$] In this case, $J$, must be $\vdash p
                    \type$, in which case there's nothing to construct.

                \item[$\textsc{Self}(\Ded)$] $J$ must be $\Gamma, J \vdash J$,
                    in which case we choose as a term the variable that
                    acompanies $J$.

                \item[$\textsc{Weak}(\Ded, \Ded')$] $J$ is \emph{anything}, with
                    one additional hypothesis in the context. We apply this
                    translation recursively in the proof of the hypothesis that
                    carries $J$ and use \emph{that} result.
                    
                \item[$\textsc{Form}_c(\Ded_0, \cdots, \Ded_{\arity c})$] $J$
                    must be of the form $\Gamma \vdash c(\cdots) \prop$, in
                    which case there's no term to build.
                    
                \item[$\textsc{Intro}^r_c(\Ded_0, \cdots, \Ded_{\arity c})$] $J$
                    must be $\Gamma \vdash C(P_i) \true$, with hypothesis of the
                    form $\Gamma \vdash P_i \true$ or $\Gamma \vdash P_i
                    \prop$. To the ones that carry the adjective $\true$, we
                    recursively apply this procedure and obtain terms $t_i$,
                    which we use to build the term $\lambda_c^r(t_i$.
                    
               \item[$\textsc{Elim}_c(\Ded, \Ded_{r_0}, \cdots, \Ded_{r_l})$]
                   $J$ must be of the form $\Gamma \vdash C \true$ with
                    hypothesis $\Gamma \vdash c(P_i) \true$ and $\Gamma,
                    P_{r_ji} \vdash C \true$, in which case we recursively apply
                    this procedure on the hypothesis and obtain terms $c_j$
                    (besides the variables $x_{r_ji}$), and construct the term
                    $\varepsilon(c_j, x_{r_ji})$.
            
            \end{description}

            If we deonte by $\term(\Ded)$ the term associated with the deduction
            $\Ded$ (if it exists), and by $\deduction(\Ded)$ the proof of its
            typing in the $\lambda$-calculus, the judgements of the
            $\lambda$-calculus will be

            \begin{align*}
                x_0 : P_0, \cdots, x_{n-1} : P_{n-1} &\vdash \term(\rho) : P \\
                x_0 : P_0, \cdots, x_{n-1} : P_{n-1} &\vdash P \type
            \end{align*}

            with the associated proofs $\deduction(\rho)$.

            Finally, we obtain the triple $(\Sigma_\lambda, J_\lambda,
            \rho_\lambda)$. Verifying that their image is what it's supposed to
            be is as simple as a matter of calculations -- which , nonetheless
            should be intuitive: as far as the signatures and judgements go,
            it's evident; as for the deductions, it's as troublesome as
            observing that the translations works rule-by-rule -- and, afterall,
            we've defined $\rho_\lambda$ so that it would work with that
            translation.

        \end{proof}

    \end{theorem}

    \citet{SorensenLectures2006}'s result is then a corollary

    \begin{corollary}[Curry-Howard]

        \foreign{Ipsis literis},

        \begin{enumerate}

            \item IF $\Gamma \vdash M : \varphi$ in
                $\lambda_\rightarrow$, then $\mathrm{rg}(\Gamma) \vdash
                \varphi$ in IPC$(\rightarrow)$.

            \item If $\Delta \vdash \varphi$ em IPC$(\rightarrow)$, then $\Gamma
                \vdash M : \varphi$ in $\lambda_\rightarrow$, for some $M$ and
                some $\Gamma$ with $\mathrm{rg}(\Gamma) = \Delta$.

        \end{enumerate}

        in our formalism: in the standard signature $\Sigma_\Prop$ of the
        implicational propositional logic and $\Sigma_\lambda$ of the
        $\lambda$-calculus,

        \begin{enumerate}

            \item If $J : \Sen_\lambda\Sigma_\lambda$ has a proof, then
                $\alpha_{\Sigma^\lambda}(J)$ has a proof, and

            \item If $J : \Sen_\Prop\Sigma_\Prop$ has a proof, then there is $J'
                : \Sen_\lambda\Sigma_\lambda$ that has a proof and
                $\alpha_{\Sigma_\lambda}(J') = J$.

        \end{enumerate}

        \begin{proof}

            If there is a proof of a judgement $J$ in the $\lambda$-calculus in
            the tradional signature $\Sigma_\lambda$, there is proof of
            $\alpha_{\Sigma_\lambda}$ in the intuitionistic propositional logic,
            as $\mathcal{T}_{\lambda\to\Prop}$ is a morphism.

            The other way, consider $J$ a judgement of the intuitionistic
            propositional logic in the standard signature $\Sigma_\Prop$, with a
            proof $\rho : \pi^\Prop_\Sigma J$. Then there is a judgement $J' :
            \Sen_\lambda\Sigma_\lambda$ and proof $\rho' :
            \pi^\lambda_{\Sigma_\lambda}J'$ with image $J$ and $\rho$ under
            $\mathcal{T}_{\lambda\to\Prop}$. In particular,
            $\alpha_{\Sigma_\lambda}J' = J$.

        \end{proof}

    \end{corollary}

    Finally, it's worth noting the following: one particular troublesome point
    about the Curry-Howard Correspondence was that, strangely, the literature
    didn't seem to agree if the theorem deserved the title of an
    \emph{isomorphism} or not. The title of \citet{SorensenLectures2006}'s is
    \textit{Lectures on the Curry-Howard Isomorphism}, even though they
    themselves note that

    \begin{quotation}

        \omissis the reader may find \inclusion{the theorem as stated} a little
        unsatisfactory. If we talk about an ``isomorphism'' then perhaps the
        statement of the proposition should have the form of an equivalence? The
        concluding sentence \inclusion{-- that ``in particular an implicational
        formula is an intuitionistic theorem if and only if it is an inhabited
        type --} is indeed of this form, but it only holds on a fairly high
        level: We must abstract from the proofs and only ask about conclusions.
        (Or, equivalently, we abstract from terms and only ask which types are
        non-empty.) To support the idea of an ``isomorphism,'' we would certainly
        prefer an exact, bijective correspondence between proofs and terms.

        Unfortunately, we cannot improve \inclusion{the statement} in this
        respect, at least not for free. While it is correct to say that
        lambda-terms are essentially annotated proofs, the problem is that some
        proofs can be annotated in more than one way. For instance, the proof
       
        \[
            \inferrule {
                \inferrule {
                    p \vdash p
                } {
                    p \vdash p \rightarrow p
                }
            } {
                \vdash p \rightarrow p \rightarrow p
            }
        \]

        can be annotated as either $\lambda x^p \lambda y^p\,x$ or $\lambda x^p
        \lambda y^p\,y$. \citep[\S~4.3]{SorensenLectures2006}

    \end{quotation}

    We may add that, further, even a single judgement of the $\lambda$-calculus
    may have more than one proof.

    \begin{mathpar}
        \inferrule {
            \inferrule {
                \inferrule {
                    \omissis
                } {
                    x : A, y : A \vdash y : A
                }
            } {
                x : A \vdash \lambda y y : A \rightarrow A
            }
        } {
            \vdash \lambda x \lambda y y : A \rightarrow (A \rightarrow A)
        } \and
        \inferrule {
            \inferrule {
                \inferrule {
                    \omissis
                } {
                    \vdash \lambda y y : A \rightarrow A
                } \\
                \inferrule {
                    \omissis
                } {
                    \vdash A \type
                }
            } {
                x : A \vdash \lambda y y : A \rightarrow A
            }
        } {
            \vdash \lambda x \lambda y y : A \rightarrow (A \rightarrow A)
        }
    \end{mathpar}

    \begin{mathpar}
        \textsc{Intro}_\rightarrow(\textsc{Intro}_\rightarrow(\omissis)) \and
        \textsc{Intro}_\rightarrow(\textsc{Weak}(\omissis))
    \end{mathpar}

    It's up to the reader to complement the proofs and verify that the trees in
    fact correspond to the written proofs.

    What our perspective may offer is that the correspondence is, in fact,
    bijective, but somehow only between the proofs in both systems. In the case
    of judgements, multiple ones may hit a particular target.

\section{Future developments}
\label{sec:further}

    There are several subject we left behind, and gaps we've left to fill. Those
    can be the object of further developments. We'll discuss them briefly.

    \subsection{Polarity}

        As you will no doudbtly have noticed, developing all the details of this
        generalized Curry-Howard Correspondence is a toilsome job. This, among
        others, was the reason for us to ommit what would be a more
        \emph{complete} formulation of propositional logic.

        You may have observed, if moved by obsessive curiosity about logical
        systems, that our formalism does not admit the usual rules for
        implication.

        \begin{mathpar}
            \inferrule {
                \Gamma \vdash A \rightarrow B \true \\
                \Gamma \vdash A \true
            } {
                \Gamma \vdash B \true
            } \and
            \inferrule {
                \Gamma, A \vdash B \true
            } {
                \Gamma \vdash A \rightarrow B \true
            }
        \end{mathpar}

        That's because our propositional logic only admits \strong{positive}
        connectives, and the implication is a \strong{negative} one. That's
        what's called the \strong{polarity} of a connective, which defines if
        the connective is determined by its introduction or elimination rules.

        As we've defined them, connectives come with a list of introduction
        rules and a single elimination rule \emph{calculated} from the former.
        Negative connectives work the other way around.

        The non-linear conjunction notoriously admits both a positive and a
        negative definition. The later is most familiar:

        \begin{mathpar}
            \inferrule {
                \Gamma \vdash A \land B \true
            } {
                \Gamma \vdash A \true
            } \and
            \inferrule {
                \Gamma \vdash A \land B \true
            } {
                \Gamma \vdash B \true
            } \and
            \inferrule {
                \Gamma \vdash A \true \\
                \Gamma \vdash B \true
            } {
                \Gamma \vdash A \land B \true
            }
        \end{mathpar}

        in which the last rule is defined from the previous by the mantra ``to
        deduce the connective, it is necessary to provide everything that from
        it may be taken''.

        The positive formulation, intead, is

        \begin{mathpar}
            \inferrule {
                \Gamma \vdash A \true \\
                \Gamma \vdash B \true
            } {
                \Gamma \vdash A \land B \true
            } \and
            \inferrule {
                \Gamma \vdash A \land B \true \\
                \Gamma, A, B \vdash C \true
            } {
                \Gamma \vdash C \true
            }
        \end{mathpar}

        where the second rule comes from the first by the mantra ``to conclude
        something from a connective, we must be able to prove it from every set
        of hypothesis that might have been used to build it''\footnotemark.

        \footnotetext{It's not a mantra that rolls out of your tongue, but it
        works.}

        Extending the formalism to admit rules of negative polarity is not hard,
        but it is laborious. We'd break the definition in two cases. For
        exemple, the headers wouldn't be simply a vector of adjectives, but
        indeed a such vector \emph{if} the connective is positive, and something
        else otherwise.

        \subsection{Universal formulation}

        Despite having formalized the Curry-Howard Correspondence, as was our
        goal, we didn't generalize satisfactorily. The step from a single
        signature to multiple ones related by (unused) morphisms is only
        obvious. A more elegant formulation would have to be more
        \emph{universal} in flavor.

        What would it look like? A candidate is to say that two p-institutions
        are in ``Curry-Howard correspondence'' if there is a morphism between
        them that's bijective in the space of all proofs.

        More analysis would be necessary, but right off the bat it seems too
        wak. Certanly there must be such morphisms that would not deserve,
        somehow, the name of Curry-Howard.

        Another possible aproach is to algorithmically build the
        $\lambda$-calculus from the propositional logic. That would hint us in
        what way to apply that procedure to more general p-institutions -- a
        sort of ``Curry-Howardization''. That construction would enojy many
        properties that, correctly abstracted, would give us something that
        looks like a universal property.

    \subsection{Other proof systems}

        Natural deduction is not alone as a manner of building proofs. Indeed,
        it is less popular than others like Hilbert-style proofs or sequents.

        All of this formalisms -- we are to believe -- form p-institutions if
        correctly parametrized. If that's the case, what kind of morphisms might
        be stablished between them? What do they tell us about the relation
        between those formalisms?

    \subsection{Other constructions}

        The language of p-institutions (paired with the powerful logic-building
        tool that are deductive systems) are generic enough to be useful not
        only in exploring the Curry-Howard Correspondence. There are several
        procedures logics are subject to that might be analysed.

        As an example, consider a tarskian relation $\succ$ over $\Judgement$.
        We define what's called it's \emph{paraconsistetization} $\succcurlyeq$.
        We call a set $\Delta \subseteq \Judgement$ \strong{inconsistent} if it
        can deduce any judgement. The paraconsistentization

        \[
            \Delta \succcurlyeq J
        \]

        stands if, and only if, some consistent\footnotemark (not inconsistent)
        subset $\Delta' \subseteq \Delta$ is such that $\Delta' \succ J$. What
        kind of p-institution construction would that look like? What universal
        properties would it enjoy?

        \footnotetext{We might as well define a ``finitarization'' of a logic in
        much the same way.}

        Furthermore, that's not the only procedure we can come up with. We could
        think of some sort of ``universe addition'', like the ones in
        Martin-Löf's Type Theory, or even something that could go
        by the name of the ``complete interior'' of a logic according to a
        certain semantic.

        The point is: a \emph{structural} exploration of logics is a fruitfull
        endeavor; and a categorical formulation -- minimalist despite expressive
        -- can be a \foreign{lingua franca} in our incursions into this subject.

    \bibliographystyle{plainnat}
    \bibliography{ref}

\end{document}